\documentclass[11pt]{amsart}
\usepackage{latexsym,amsfonts,amssymb}
\newtheorem{theorem}{Theorem}[section]
\newtheorem{lemma}[theorem]{Lemma}
\newtheorem{corollary}[theorem]{Corollary}

\theoremstyle{definition}

\theoremstyle{remark}

\numberwithin{equation}{section}

\newcommand{\ab}{_{\alpha \bar{\beta}}}
\begin{document}
\noindent
\vspace{0.5in}

\title{Non-ancient solution of the Ricci flow}

\author{Qihua Ruan$^{1}$}
\author{Zhihua Chen$^{2}$}
\address{$^{1,2}$Department of Applied Mathematics,Tongji University,
Shanghai,200092 P.R.China.}
\address{$^{1}$Department of Mathematics, Putian College, Putian, Fujian 352100 P.R.China.}
 \email{ruanqihua@163.com}
\thanks{Project supported by NFSC No.
10271089 ;and STC of Shanghai No.03JC1402. }


\begin{abstract}  For any complete noncompact K$\ddot{a}$hler manifold with nonnegative and
bounded holomorphic bisectional curvature,we provide the necessary
and sufficient  condition for   non-ancient solution to the Ricci
flow in this paper.
\end{abstract}

\maketitle

\section{Introduction}
Let M be any complex n-dimensional complete noncompact
K$\ddot{a}$hler manifold with nonnegative and bounded holomorphic
bisectional curvature, in this paper we always assume M satisfies
this condition.The Ricci flow
$$\frac{\partial}{\partial t}g_{\alpha \bar{\beta}}=-R_{\alpha \bar{\beta}},\ \ \
g_{\alpha \bar{\beta}}(x,0)=g_{\alpha \bar{\beta}}(x)
\leqno(1.1)$$ evolves K$\ddot{a}$hler metric $g\ab$ on M by its
Ricci tensor $R\ab$. It was first introduced by Hamilton
(\cite{H1}),in (\cite{H2}) he divided a solution to the Ricci flow
into three types, the first two types was called ancient solution,
we call the last type by non-ancient solution,which is said that
the Ricci flow has long time existence an $0\leq R(x,t)\leq
\frac{C}{1+t}$,where R(x,t)denotes the scalar curvature at time t
and C is a positive constant.Chen-Zhu (\cite{CZ1}) proved
non-ancient solution exists if M is a complete noncompact complex
two-dimensional K$\ddot{a}$hler manifold with positive and bounded
holomorphic bisectional curvature,its geodesic balls have
Euclidean volume growth and its scalar curvature decays to zero at
infinity in the average sense. They used this result to study the
uniformization conjecture by Yau(\cite{Y}), which says that a
complete noncompact K$\ddot{a}$hler manifold of positive
holomorphic bisectional curvature is biholomorphic to a complex
Euclidean space. They partially confirmed this conjecture in the
case of n=2. Recently Ni and Tam (\cite{NT}) studied the Ricci
flow by solving the Poincare-Lelong equation, and got some nice
results.

In this paper, we combined with Ni and Tam's results  to obtain
the necessary and sufficient condition for the non-ancient
solution to the Ricci flow. That is :
\begin{theorem} Let M be above assumption ,then the Ricci flow
(1.1) has non-ancient solution if and only if
$$\int^{r}_{0}sk(x,s)ds\leq C log (2+r)\leqno(1.2)$$
for some constants $C >0,\, \forall x \in M , r\geq 0$, where
 $$k(x,\,s)=\frac{1}{V(x,s)}\int_{B(x,\,s)}R(x)dx$$
 V(x,s) is the volume of the geodesic ball B(x,\,s) centered at $x\in
 M$ with radius $s$,R(x) is the scalar curvature of M.
\end{theorem}

 \vspace{0.3in}
\section{The necessary condition for the solution}
For the convenience, we  suppose $C>0 $ be a various constant
depending only on $n$ in the following.
\begin{lemma}(\cite{S})Let M be above assumption, then there
exists a constant C such that for $\forall x,y \in M,$
$$\frac{d(x,\,y)^{2}}{C V(x,\,d(x,\,y))}\leq G(x,\,y)\leq \frac{C d(x,\,y)^{2}}{V(x,\,d(x,\,y))} \leqno(2.1)$$
$$\mid\nabla G(x,\,y)\mid\leq \frac{C d(x,\,y)}{V(x,\,d(x,\,y))} \leqno(2.2)$$
where $d(x,\,y)$ denotes the distance between $x$ and $y$ with
respect to $g\ab(x)$, $\nabla$ is the covariant derivative with
respect to $g\ab(x)$, $G(x,\,y)$ is the positive Green's function
on $M$ with respect to $g\ab(x)$.
\end{lemma}
\begin{theorem} Let $M$ be above assumption, and the necessary condition holds, then there exists a
constant $C>0$ such that for $\forall x \in M ,r>0,$ we have
$$\int_{B(x,\,r)}\frac{R(y)d(x,\,y)^{2}}{V(x,\,d(x,\,y))}dy \leq C\, log(2+r) $$
\end{theorem}
\begin{proof}
Let $g\ab(x,\,t)$ be the solution of the Ricci flow (1.1) with
$g\ab(x)$ as the initial metric. From the necessary condition, we
know that the solution exists for all the times and satisfies
$$0\leq R(x,\,t)\leq \frac{C}{1+t} \,\, \,\,on \,\,M\times [0,\,+\infty) \leqno(2.3)$$
Let $F(x,\,t)= log \frac{det(g\ab(x,\,t))}{det(g\ab(x,0))}$, since
$-\partial_{\alpha}\bar{\partial_{\beta}}log
\frac{det(g_{\delta\bar{\gamma}}(x,\,t))}{det(g_{\delta\bar{\gamma}}(x,\,0))}
=R\ab(x,\,t)-R(x,\,0)$, after taking trace with $g\ab(x,\,0)$, we
get
$$R(x,\,0)= \triangle F(x,\,t)+ g^{\alpha\bar{\beta}}(x,\,0)R\ab(x,\,t) \leqno(2.4)$$
where $\triangle$ is the Laplace operator of the metric
$g\ab(x,\,0)$.

For any fixed $x_{0}\in M $ and any $\alpha>0$, we denote
$$\Omega_{\alpha}=\{x\in M | G(x_{0},\,x)\geq \alpha \}$$
combining (2.1) with Shi's technique on $M\times C^{2}$(\cite{S}),
it is not hard to see that there exist a number $d(\alpha)\geq 1$
such that
$$\frac{d(\alpha)^{2}}{V(x_{0},\,d(\alpha))}=\alpha \leqno(2.5)$$
and a constant $C>0$ such that
$$B(x_{0},\,C^{-1} d(\alpha))\subset \Omega_{\alpha}\subset B(x_{0},\,C
d(\alpha))\leqno(2.6)$$ Recall that $F(x,\,t)$ evolves by
$$\frac{\partial F(x,\,t)}{\partial t}=-R(x,\,t)$$
Combining with (2.3), we obtain
$$0\geq F(x,\,t)\geq -C\,log (1+t)\leqno(2.8)$$
Multiplying (2.7) by $G(x_{0},\,x)-\alpha$ and integrating over
$\Omega_{\alpha}$, we have

 $$  \int_{\Omega_{\alpha}}R(x,\,0)(G(x_{0},\,x)-\alpha)dx \,\,\,\,\,\,\,\,\,\,\,\,\,\,\,\,\,\,\,\,\,\,\,\,
 \,\,\,\,\,\,\,\,\,\,\,\,\,\,\,\,\,\,\,\,\,\,\,\,\,\,\,\,\,\,\,\,\,\,\,\,\,\,\,\,
 \,\,\,\,\,\,\,\,\,\,\,\,\,\,\,\,\,\,\,\,\,\,\,\,\,\,\,\,\,\,\,
  \,\,\,\,\,\,\,\,\,\,\,\,\,\,\,\,\leqno(2.9)$$
 $$=\int_{\Omega_{\alpha}}(\triangle
F(x,\,t))(G(x_{0},\,x)-\alpha)dx+\int_{\Omega_{\alpha}}g^{\alpha\bar{\beta}}(x,\,0)R\ab(x,\,t)(G(x_{0},\,x)-\alpha)dx$$
$$=-\int_{\partial \Omega_{\alpha}}F(x,\,t)\frac{\partial G}{\partial \nu}d\sigma - F(x_{0},\,t) +
\int_{\Omega_{\alpha}}g^{\alpha\bar{\beta}}(x,\,0)R\ab(x,\,t)(G(x_{0},\,x)-\alpha)dx$$
$$\leq C (1+\int_{\partial \Omega_{\alpha}}\frac{\partial G}{\partial \nu}d\sigma) log (1+t)
+
\int_{\Omega_{\alpha}}g^{\alpha\bar{\beta}}(x,\,0)R\ab(x,\,t)G(x_{0},\,x)dx\,\,\,\,\,\,\,\,\,\,\,\,\,
\,\,\,\,\,\,\,\,\,\,\,\,$$ Here we use (2.8) and denote $\nu$ by
the outer unit normal of $\partial\Omega_{\alpha}$.
\\Integrating (2.9) from $\alpha$ to $2\alpha$ and multiplying by
$\frac{1}{\alpha}$, we get
$$\int_{\Omega_{2\alpha}}R(x,\,0)(G(x_{0},\,x)-2\alpha)dx\,\,\,\,\,\,\,\,\,\,\,\,\,\,\,\,\,\,\,\,\,\,\,\,
 \,\,\,\,\,\,\,\,\,\,\,\,\,\,\,\,\,\,\,\,\,\,\,\,\,\,\,\,\,\,\,\,\,\,\,\,\,\,\,\,
 \,\,\,\,\,\,\,\,\,\,\,\,\,\,\,\,\,\,\,\,\,\,\,\,\,\,\,\,\,\,\,
  \,\,\,\,\,\,\,\,\,\,\,\,\,\,\,\,\leqno(2.10)$$
  $$\leq C (1+\frac{1}{\alpha}\int^{2\alpha}_{\alpha}\int_{\partial \Omega_{\gamma}}\frac{\partial G}{\partial \nu}d\sigma d \gamma) log (1+t)
+
\int_{\Omega_{\alpha}}g^{\alpha\bar{\beta}}(x,\,0)R\ab(x,\,t)G(x_{0},\,x)dx$$
It is easy to see that
$$\frac{1}{2}\int_{\Omega_{4\alpha}}R(x,\,0)G(x_{0},\,x)dx \leq\int_{\Omega_{2\alpha}}R(x,\,0)(G(x_{0},\,x)-2\alpha)dx\leqno(2.11)$$
Using Shi's technique on $M \times C^{2}$ and combining (2.1) with
(2.5), we obtain that there exists a constant $C>0$ such that
$$C^{-1}d(\alpha)\leq d(\gamma)\leq C d(\alpha), \alpha\leq\gamma\leq2\alpha.\leqno(2.12)$$
From the coarea formula, we have that
$$d\sigma d\gamma = \frac{\partial G(x_{0},\,x)}{\partial \nu}d \sigma d\nu
=|\frac{\partial G(x_{0},\,x)}{\partial \nu}|d \sigma
|d\nu|=|\frac{\partial G(x_{0},\,x)}{\partial
\nu}|dx\leqno(2.13)$$ Combining (2.2),(2.5),(2.6),(2.12),(2.13)
and the standard volume comparison, we have that
$$\frac{1}{\alpha}\int^{2\alpha}_{\alpha}\int_{\partial \Omega_{\gamma}}\frac{\partial G(x_{0},\,x)}{\partial \nu}d\sigma d \gamma
\leq \frac{1}{\alpha}\int^{2\alpha}_{\alpha}\int_{\partial
\Omega_{\gamma}}|\frac{\partial G(x_{0},\,x)}{\partial
\nu}|^{2}dx\leqno(2.14)$$
$$\leq \frac{C}{\alpha}\int ^{2\alpha}_{\alpha}\int_{\partial \Omega_{\gamma}}(\frac{d(\gamma)}{V(x_{0},\,d(\gamma))})^{2}dx
\leq\frac{C d(\alpha)^{2}}{\alpha
(V(x_{0},\,C^{-1}d(\alpha)))^{2}}\int _{\Omega_{\alpha}\setminus
\Omega_{ 2\alpha}}dx$$
$$\leq\frac{C d(\alpha)^{2}}{\alpha
(V(x_{0},\,C^{-1}d(\alpha)))^{2}} V(x_{0},\,Cd(\alpha))\leq C
\,\,\,\,\,\,\,\,\,\,\,\,\,\,\,\,\,\,\,\,\,\,\,\,\,\,\,\,\,\,\,\,\,\,\,\,\,\,\,\,\,\,\,\,\,\,\,\,\,\,\,\,\,\,\,\,\,
\,\,\,\,\,$$ Substituting (2.14) and (2.11) into (2.10), we get
that
$$\int_{\Omega_{4\alpha}}R(x,\,0)G(x_{0},\,x)dx \leq C log (1+t)+\int_{\Omega_{\alpha}}g^{\alpha\bar{\beta}}(x,\,0)R\ab(x,\,t)G(x_{0},\,x)dx$$
Integrating this inequality from $0$ to $t$, we see that for any
$t>0$,$$\int_{\Omega_{4\alpha}}R(x,\,0)G(x_{0},\,x)dx \leq C log
(1+t)+\frac{1}{t}\int_{\Omega_{\alpha}}g^{\alpha\bar{\beta}}(x,\,0)(g\ab(x,\,0)-g\ab(x,\,t))G(x_{0},\,x)dx$$
$$\leq C log(1+t) +
\frac{C}{t}\int_{\Omega_{\alpha}}G(x_{0},\,x)dx\,\,\,\,\,\,\,\,\,\,
$$
Finally,combining this inequality with (2.1) and (2.6), we have
that there exists a constant $C>0$ such that for $\forall x_{0}\in
M, t>0 $ and $r>0$
$$\int_{B(x_{0},\,r)}R(x,\,0)\frac{d(x_{0},\,x)}{V(x_{0},\,d(x_{0},\,x))}dx\leq C (log(1+t)+\frac{r^{2}}{t})$$
Choose $t=r^{2}$, we can prove Theorem 2.2 .
\end{proof}
\vspace{0.3in}
\begin{lemma}(\cite{NST}): Let $M$ be above assumption,then the
Poisson equation  $\triangle u(x) =R(x)$ has a solution with
$\sup\limits_{x\in B(x_{0},\,r)}|u(x)|\leq C log(2+r)$ for a
constant $C>0, \forall r>0$, if and only if $$\int
^{t}_{0}sk(x,\,s)ds \leq C' log(2+t)$$ for a constant $ C'>0,
\forall t\geq \frac{1}{5}r.$
\end {lemma}

In the following, we prove the necessary part.
\begin{proof}
From Lemma 2.3, it is  sufficient to prove that the Poisson
equation $\triangle u(x)=R(x)$ has a solution with
$\sup\limits_{x\in B(x_{0},\,r)}|u(x)|\leq C log(2+r)$ . In
fact,if $t\geq \frac{1}{5}r$, then from Lemma 2.3 we can conclude
the necessary condition.  if $t\leq \frac{1}{5}r$, since the
scalar curvature is bounded, then $ \int ^{t}_{0}sk(x,\,s)ds \leq
C$ , so there always exists a constant $C'>0$ such that $\int
^{t}_{0}sk(x,\,s)ds \leq C' log(2+t)$ for all $t\geq 0$.

To solve the poisson equation, we first construct a family of
approximate solution $u_{r}$ as follows.

For a fixed $x_{0} \in M$ and $\forall r > 0,$ define $u_{r}(x)$
on $B(x_{0},\,r)$ by
$$u_{r}(x) = \int_{B(x_{0},\,r)}(G(x_{0},\,y) - G(x,\,y))R(y)d(y)$$
It is clear that $u_{r}(x_{0}) = 0,$ and $\triangle u_{r}(x) =
R(x)$ on $B(x_{0},\,r).$ For $x \in B(x_{0},\,\frac{r}{4}),$ we
write
$$u_{r}(x) = (\int_{B(x_{0},\,r)\setminus
B(x_{0},\,4d(x_{0},\,x))} +
\int_{B(x_{0},\,4d(x_{0},\,x))})(G(x_{0},\,y) - G(x,\,y))R(y)dy =
I_{1} + I_{2}$$ From Theorem2.2, we see that
$$\mid I_{2} \mid \leq C\log(2+d(x,\,x_{0})) \,\,\,\,\mbox{on}\,\,\,\,B(x_{0},\,\frac{r}{4})\leqno(2.15)$$
From(2.2), we have that for $y \in B(x_{0},\,r) \setminus
4d(x_{0},\,x),$
$$\mid G(x_{0},\,y) - G(x,\,y) \mid \leq d(x,\,x_{0})\sup\limits_{z \in B(x_{0},\,d(x,\,x_{0}))}
\mid
\nabla_{z}G(z,\,y)\mid\,\,\,\,\,\,\,\,\,\,\,\,\,\,\,\,\,\,\,\,\,\,\,\,\,\,\,\,\,\,\,\,\,
\,\,\,\,\,\,\,\,\,\,\,\,\,\,\,\,\,\,\,\,\,\,\,\,\,\,\,\,$$

$$\leq Cd(x,\,x_{0})\sup\limits_{z \in B(x_{0},\,d(x,\,x_{0}))}\frac{d(z,\,y)}{V(z,\,d(z,\,y))}$$

$$\leq C\frac{d(x,\,x_{0})d(x_{0},\,y)}{V(x_{0},\,d(x_{0},\,y))}
\,\,\,\,\,\,\,\,\,\,\,\,\,\,\,\,\,\,\,\,\,\,\,\,\,\,\,\,\,\,\,\,\,\,\,\,\,\,\,\,\,\,\,\,\,\,\,
\,\,\,$$

Combining this with Theorem2.2, we have that:

$$\mid I_{1}\mid \leq Cd(x,\,x_{0})\int_{B(x_{0},\,r) \setminus 4d(x_{0},\,x)}\frac{R(y)d(x_{0},\,y)}{V(x_{0},\,d(x_{0},\,y))}dy
\,\,\,\,\,\,\,\,\,\,\,\,\,\,\,\,\,\,\,\,\,\,\,\,\,\,\,\,\,\,\,\,\,\,\,\,\,\,\,\,\,\,\,\,\,\,\,\,\,\,\,\,\,\,\,\,\,\,\,\,\,
\,\leqno(2.16)$$

$$\leq Cd(x_{0},\,x)\sum_{k=2}^{\infty}\frac{1}{2^{k}d(x,\,x_{0})}\int_{B(x_{0},\,2^{k+1}d(x,\,x_{0}))\setminus
B(x_{0},\,2^{k}d(x,\,x_{0}))}\frac{R(y)d(x_{0},\,y)^{2}}{V(x_{0},\,d(x_{0},\,y))}dy$$

$$\leq C\sum_{k=2}^{\infty}\frac{1}{2^{k}}\log(2+2^{k+1}d(x,\,x_{0}))\,\,\,\,\,\,\,\,\,\,\,\,\,\,\,\,\,\,\,\,
\,\,\,\,\,\,\,\,\,\,\,\,\,\,\,\,\,\,\,\,\,\,\,\,\,\,\,\,\,\,\,\,\,\,\,\,\,\,\,\,\,\,\,\,\,\,\,\,\,\,\,\,\,
\,\,\,\,\,\,\,\,\,\,\,\,\,\,\,\,\,\,\,\,\,\,\,\,\,\,\,\,\,\,\,\,\,\,\,
\,\,\,\,\,\,\,\,\,\,$$

$$\leq C\log(2+d(x,\,x_{0}))\,\,\,\,\,\,\,\,\,\,\,\,\,\,\,\,\,\,\,\,\,\,\,\,\,\,\,\,\,\,\,\,\,\,\,\,\,\,\,\,\,\,
\,\,\,\,\,\,\,\,\,\,\,\,\,\,\,\,\,\,\,\,\,\,\,\,\,\,\,\,\,\,\,\,\,\,\,\,\,\,\,\,\,\,\,\,\,\,\,\,\,\,\,\,\,\,\,\,\,
\,\,\,\,\,\,\,\,\,\,\,\,\,\,\,\,\,\,\,\,\,\,\,\,\,\,\,\,\,\,\,\,\,\,\,\,\,\,\,\,\,\,\,\,\,\,\,\,\,$$

Hence, from (2.15) and (2.16), we deduce that
$$\mid u_{r}(x)\mid \leq C \log (2+d(x,\,x_{0}))\,\,\,\,\,\mbox{for}\,\,\,\,\,r \geq 4d(x_{0},\,x).$$

Therefore, it follows from the Schauder theory of elliptic
equations that there exists a sequence of $r_{j}\rightarrow
+\infty$ such that $u_{r_{j}}(x)$ converges uniformly on compact
subset of $M$ to a smooth function $u$ satisfying
$$
\left\{\begin{array}{ll} u(x_{0})=0\,\,and\,\,\,\,\triangle u = R\,\,\,\,on M \\
\mid u(x) \mid \leq C \log(2+d(x,\,x_{0})),\,\,for \,\,x \in M
\end{array}\right.$$

Thus we proved the necessary part of Theorem1.1.
\end{proof}

 \vspace{0.3in}
\section{The sufficient condition of the solution}

\begin{lemma}(\cite{NST}):\,\,Let M be above assumption, if there
exists a constant $C > 0,\,\varepsilon > 0,$ such that
$$ k(x_{0},\,r) \leq \frac{C}{(1+r)^{1+\varepsilon}}, \,\,\,\,\,\,\mbox{for}\,\,\,\,x_{0}\in M,\,\,\,\,\forall r \geq 0$$
\end{lemma}
then the Ricci flow (1.1) has long time existence, $R(x,\,t)$ is
nonnegative and bounded.

\begin{lemma}(\cite{H2}):\,\,Let M be above assumption, then there
exists a constant $C > 0,$ such that for all $(x_{0},\,t) \in M
\times [0,\,+\infty)$
$$ -F(x_{0},\,t) \leq C[(1+\frac{t(1-m(t))}{R^{2}})\int^{R}_{0}sk(x_{0},\,s)ds-
\frac{tm(t)(1 - m(t))}{R^{2}}]$$ where $m(t) = \inf_{x \in
M}F(x,\,t)$
\end{lemma}

Now we can use Lemma3.1 and Lemma3.2 to prove the sufficient part.
\begin{proof}
From the sufficient condition of Theorem1.1, we have
$$\int_{0}^{r}sk(x,\,s)ds \leq C \log(2+r)$$
On the other hand,
$$\int_{0}^{r}sk(x,\,s)ds \geq
\int_{0}^{r}s\frac{1}{V(x,\,s)}\int_{B(x,\,s)}R(y)dyds$$
$$\,\,\,\,\,\,\,\,\,\,\,\,\,\,\,\,\,\,\,\,\,\,\,\,\,\,\,\,\,\,\,\,\,\,\,\,\,\,\,\,\geq
\frac{1}{V(x,\,r)}\int_{\frac{r}{2}}^{r}s\int_{B(x,\,s)}R(y)dyds$$
$$\,\,\,\,\,\,\,\,\,\,\,\,\,\,\,\,\,\,\,\,\,\,\,\geq
\frac{\frac{3}{4}r^{2}}{V(x,\,r)}\int_{B(x,\,\frac{r}{2})}R(y)dy$$

From the standard volume comparison, we have that there exists a
constant $C$ and $\varepsilon > 0$ such that
$$k(x_{0},\,r) \leq  \frac{C}{(1+r)^{1+\varepsilon}},\,\,\,\,\mbox{for}\,\,\,\,\,x_{0}\in M,\,\,\,\,\forall r \geq 0.$$
Thus form Lemma3.1, we know that the Ricci flow has long time
solution.

From Lemma3.2, we have
$$-m(t) \leq C[(1+\frac{t(1-m(t))}{R^{2}})\int^{2R}_{0}sk(x_{0},\,s)ds-
\frac{tm(t)(1 - m(t))}{R^{2}}]$$ Let $R^{2} = 2Ct(1 - m(t)),$ from
the sufficient condition, we have
$$-m(t) \leq C\log(2+2R)\leqno(3.1)$$

Since $\frac{\partial F(x,\,t)}{\partial t} =
-R(x,\,t),\,R(x,\,t)$ is bounded and $F(x,\,0) = 0,$ then
$-F(x,\,t) = \int_{0}^{t}R(x,\,s)ds \leq Ct,\,$ for all $x \in
M,\,t \geq 0;$ So
$$-m(t)\leq Ct \leqno(3.2)$$
from(3.2), we have
$$R^{2} = 2Ct(1 - m(t)) \leq Ct^{2}$$

$$R \leq Ct \leqno(3.3)$$

Substituting(3.3) into (3.1), we get
$$-m(t) \leq C\log(2+t)\leqno(3.4)$$
By Li-Yau-Hamilton type inequality[8]: $\frac{\partial R}{\partial
t}+\frac{R}{t} \geq 0,$ we know $\frac{\partial
tR(x,\,t)}{\partial t} \geq 0,$ then for all $t \geq T,$ we have
$$TR(x,\,T) \leq tR(x,\,t)$$
$$TR(x,\,T)\frac{1}{t} \leq R(x,\,t)$$
$$\int_{T}^{t}TR(x,\,T)\frac{1}{s} ds \leq \int_{T}^{t}R(x,\,s)ds \leq -F(x,\,t) \leq C\log(2+t)$$
$$TR(x,\,T)(\log t - \log T) \leq C \log(2+t)$$
Let $t\rightarrow +\infty,$ we see that $R(x,\,T) \leq
\frac{C}{1+T},$ for some constant $C >0,$ all $T \geq 0.$

So the Ricci flow has non-ancient solution.
\end{proof}
\begin{corollary} Let
M be above assumption, if $\int_{0}^{r}sk(x,\,s)ds \leq C
\log(2+r)$ for all $x \in M,\,r \geq 0,$ and $V(x,\,r)\geq
Cr^{2n}$ for some constants $C > 0,\,r \geq 0,$ then M is
diffeomorphic to $R^{2n}$ in case $n \geq 3.$
\end{corollary}

\begin{proof}
From Theorem1.1, we know that the Ricci flow has long time
existence and $0 \leq R(r,\,t) \leq \frac{C}{1+t},$Let $t
\rightarrow +\infty,$ we have $R(x,\,t)\rightarrow +\infty,$it
means that the Ricci flow will improves the injectivity radius to
$\infty$ along the flow. The rest argument is the same as in
section3 of \cite{Y}.

\end{proof}

\bibliographystyle{amsplain}

\end{document}